\documentclass[11pt]{article}%
\usepackage[english]{babel}
\usepackage{indentfirst}
\usepackage{latexsym}
\usepackage{amsbsy,amsmath,amsfonts,amssymb,amsthm, graphics}
\usepackage{times}
\usepackage{graphicx}
\usepackage[english]{babel}
\usepackage[latin1]{inputenc}
\usepackage[T1]{fontenc}
\usepackage[usenames]{color}
\usepackage[colorlinks=true, bookmarksnumbered=true, bookmarksopen=true,
bookmarksopenlevel=3, pdfstartview=FitH, linkcolor=cyan, pdfmenubar=true,
pdftoolbar=true, bookmarks=true,citecolor=cyan, urlcolor=magenta,
filecolor=cyan,menucolor=black,plainpages=false,pdfpagelabels]{hyperref}
\usepackage[lmargin=3.2cm,tmargin=3cm,bmargin=3cm,rmargin=3.2cm]{geometry}
\usepackage{amsmath}
\usepackage{amsfonts}
\usepackage{amssymb}%
\numberwithin{equation}{section}
\newtheorem{tma}{Theorem}[section]

\newtheorem{obs}[tma]{Remark}
\newtheorem{lem}[tma]{Lemma}

\def\fin { \vskip 0pt \hfill \hbox{\vrule height 5pt width 5pt depth 0pt} \vskip 12pt}
\begin{document}

\title{Global smoothness for a 1D supercritical transport model with nonlocal velocity}
\author{{Lucas C. F. Ferreira$^{1}$}{\thanks{{Corresponding author.}\newline{E-mail
addresses: lcff@ime.unicamp.br (L.C.F. Ferreira), valtermoitinho@live.com
(V.V.C Moitinho)}\newline{LCF Ferreira was supported by FAPESP and CNPq,
Brazil.\newline V.V.C Moitinho was supported by CAPES and CNPq, Brazil.}}},
\ {Valter V. C. Moitinho$^{2}$}\\\\{\small $^{1,2}$ State University of Campinas (Unicamp), IMECC-Department of
Mathematics,} \\{\small Rua S\'{e}rgio Buarque de Holanda, 651, CEP 13083-859, Campinas, SP,
Brazil.}}
\date{}
\maketitle

\begin{abstract}
We are concerned with a nonlocal transport 1D-model with supercritical
dissipation $\gamma\in(0,1)$ in which the velocity is coupled via the Hilbert
transform, namely the so-called CCF model. This model arises as a lower
dimensional model for the famous 2D dissipative quasi-geostrophic equation and
in connection with vortex-sheet problems. It is known that its solutions can
blow up in finite time when $\gamma\in(0,1/2)$. On the other hand, as stated
by Kiselev (2010), in the supercritical subrange $\gamma\in\lbrack1/2,1)$ it
is an open problem to know whether its solutions are globally regular. We show
global existence of non-negative $H^{3/2}$-strong solutions in a supercritical
subrange (close to 1) that depends on the initial data norm. Then, for each
arbitrary smooth non-negative initial data, the model has a unique global
smooth solution provided that $\gamma\in\lbrack\gamma_{1},1)$ where
$\gamma_{1}$ depends on the $H^{3/2}$-initial data norm. Our approach is
inspired on that of Coti Zelati and Vicol (IUMJ, 2016).

{\small \bigskip\noindent\textbf{AMS MSC:} 35Q35; 35B65; 76D03}

{\small \medskip\noindent\textbf{Keywords:} 1D transport model; Nonlocal
velocity; Hilbert transform; Global regularity; Supercritical dissipation}

\end{abstract}

\renewcommand{\abstractname}{Abstract}

\section{Introduction}

We consider the initial value problem for the $1$D transport equation with
nonlocal velocity%

\begin{equation}
\left\{
\begin{split}
&  \partial_{t}\theta-\mathcal{H}\theta\theta_{x}+\Lambda^{\gamma}%
\theta=0\text{ \ in }\ \mathbb{T}\times(0,\infty)\\
&  \theta(x,0)=\theta_{0}(x)\text{ \ in }\ \mathbb{T},
\end{split}
\right.  \label{1Dmodel}%
\end{equation}
where $0<\gamma\leq2$, $\mathbb{T}$ is the $1$D torus, $\Lambda=(-\Delta
)^{1/2}$ and $\mathcal{H}$ denotes the Hilbert transform.

In the literature, the equation (\ref{1Dmodel}) arises as an one-dimensional
model for the famous 2D dissipative quasi-geostrophic equation (2DQG) and in
connection with vortex sheet evolution (see e.g. \cite{Baker-Morlet-1,
Cordoba}). For more details and results about the 2DQG see
\cite{Constantin1,Cordoba2,Vicol,Kiselev} and references therein.

In view of the transport structure of (\ref{1Dmodel}), any sufficiently
regular solution satisfies the following maximum principle
\[
\left\Vert \theta(\cdot,t)\right\Vert _{L^{\infty}}\leq\left\Vert \theta
_{0}\right\Vert _{L^{\infty}(\mathbb{T})}.
\]
For non-negative initial data $\theta_{0}$, one can show that the solution
$\theta$ is also non-negative. Moreover, the condition $\theta_{0}\geq0$ is
required to obtain
\begin{equation}
\left\Vert \theta(\cdot,t)\right\Vert _{L^{2}(\mathbb{T})}\leq\left\Vert
\theta_{0}\right\Vert _{L^{2}(\mathbb{T})} \label{max-L2}%
\end{equation}
for all $t\geq0$, by means of energy methods (see \cite{Cordoba}).

The IVP (\ref{1Dmodel}) has three basic cases: subcritical $1<\gamma\leq2$,
critical $\gamma=1$ and supercritical $0<\gamma<1$. The global smoothness
problem in the critical and subcritical cases have already been solved. In
\cite{Cordoba}, C\'{o}rdoba, C\'{o}rdoba and Fontelos proved global regularity
for non-negative initial data $\theta_{0}\in H^{2}$ for $1<\gamma\leq2.$ In
the critical case, they obtained global $H^{1}$-solutions by assuming
smallness condition on $L^{\infty}$-norm of the non-negative initial data.
Also, in the inviscid case of (\ref{1Dmodel}), i.e., without the viscous term
$\Lambda^{\gamma}\theta,$ they showed blow-up of solutions for compactly
supported, even and non-negative $C^{1+\varepsilon}(\mathbb{R})$-initial data
such that $\max_{x\in\mathbb{R}}\theta_{0}(x)=\theta_{0}(0)=1$. In
\cite{Dong}, Dong showed global well-posedness of (\ref{1Dmodel}) for
arbitrary initial data in $H^{s_{0}}$ where $s_{0}=\max\{\frac{3}{2}%
-\gamma,0\}$ and $1\leq\gamma\leq2$ (critical and subcritical cases). In the
supercritical case, he assumed a smallness condition on the initial data.

The global regularity problem for solutions of (\ref{1Dmodel}) in the
supercritical case remains an open problem. In the part $0<\gamma<\frac{1}{2}$
of the supercritical range, Li and Rodrigo \cite{Li} proved blow-up of
solutions in finite time for non-positive, smooth, even and compactly
supported initial data satisfying $\theta_{0}(0)=0$ and a suitable weighted
integrability condition. In \cite{Kiselev2}, still in the same range, Kiselev
showed blow-up of solutions in finite time for even, positive, bounded and
smooth initial data $\theta_{0}$ satisfying $\max_{x\in\mathbb{R}}\theta
_{0}(x)=\theta_{0}(0)$ and suitable integrability conditions.

In the range $\frac{1}{2}\leq\gamma<1$, to the best of our knowledge, the
formation of singularity in finite time or global smoothness is an open
problem (stated by \cite[p.~251]{Kiselev2}), even for sign restriction on the
initial data, i.e., $\theta_{0}\geq0$ or $\theta_{0}\leq0$. In \cite{Do}, for
$0<\gamma<1,$ Do obtained eventual regularization of solutions for
non-negative initial data. He also obtained global regularity for a modified
1D model with $\Lambda^{\gamma}\theta$ replaced by
\[
\mathcal{L}(\theta)=-\frac{\Lambda}{\log(1-\Delta)}\theta,
\]
which can be understood as supercritical dissipation in a $\log$-sense. In
\cite{Silvestre-Vicol}, Silvestre and Vicol provided four essentially
different proofs of blow-up of solutions in the inviscid case. Moreover, they
conjectured that solutions obtained as vanishing viscosity approximations
could be bounded in $C^{1/2},$ for all $t>0,$ which would possibly yield
H\"{o}lder regularization effects for the case $1/2\leq\gamma<1$ and then
would solve the global regularity conjecture in \cite[p.~251]{Kiselev2} (see
Conjectures 7.1 and 7.2 in \cite{Silvestre-Vicol}). In \cite{Bae}, Bae,
Granero-Belinch\'{o}n and Lazar considered the inviscid case and developed a
theory of global weak super solutions for (\ref{1Dmodel}) with non-negative
data $\theta_{0}\in L^{1}\cap L^{2}.$

In this paper we focus on supercritical values of $\gamma$ contained in the
range $\frac{1}{2}\leq\gamma<1$ (in fact, close to $1$) and prove existence of
global classical solutions for (\ref{1Dmodel}). More precisely, we show
existence of $H^{\frac{3}{2}}$-strong solution for arbitrary non-negative
initial data $\theta_{0}\in H^{\frac{3}{2}}$ and $\gamma_{1}\leq\gamma<1$,
where $\gamma_{1}$ depends on the $H^{\frac{3}{2}}$-norm of $\theta_{0}.$
Indeed, due to standard regularization of $H^{\frac{3}{2}}$-strong solutions,
our solutions are $C^{\infty}$-smooth for $t>0.$ For that matter, first we
obtain an eventual regularization result with an explicit control on the
regularization time $T^{\ast}$, namely
\[
T^{\ast}=C\alpha^{\frac{1}{1-\gamma}}\left\Vert \theta_{0}\right\Vert
_{L^{\infty}(\mathbb{T})}^{\frac{\gamma}{1-\gamma}}.
\]
Afterwards, we provide an explicit estimate for the local existence time
$T_{1}$ of $H^{\frac{3}{2}}$-solutions and then compare it with $T^{\ast}.$
More precisely, we obtain the explicit control
\[
T_{1}=C\left(  \left\Vert \theta_{0}\right\Vert _{L^{2}}^{\frac{2\gamma
(9+2\gamma)}{3(9+4\gamma)}}\left\Vert \theta_{0}\right\Vert _{\dot{H}%
^{\frac{3}{2}}}^{2-\frac{4\gamma(6+\gamma)}{3(9+4\gamma)}}\right)  ^{-1}.
\]
Let us remark that in view of (\ref{max-L2}) the non-negative condition is
relatively common for (\ref{1Dmodel}) and have been assumed in several works
(see, e.g., \cite{Bae, Bae2, Cordoba, Dong, Kiselev2}). Also, it is worth
mentioning that, in consonance with the conjectures in \cite{Silvestre-Vicol},
we obtain boundedness of solutions in $C^{\alpha}$ for $\alpha>1-\gamma$.

Coti Zelati and Vicol \cite{Vicol} showed existence and uniqueness of global
$H^{2}$-strong solutions (classical) for the 2DQG with supercritical values of
$\gamma$. Our approach follows the spirit of \cite{Vicol}, mainly with respect
to the control on the regularization time of our $H^{\frac{3}{2}}$-solutions.
However, for the control on the local existence time $T_{1}$ we need to
proceed in a different way by employing commutator estimates and suitable
energy estimates for the $H^{\frac{3}{2}}$-solutions (see estimates
(\ref{local1}) to (\ref{local7})).

Our main result reads as follows.

\begin{tma}
\label{global} Let $\theta_{0}\in H^{\frac{3}{2}}(\mathbb{T})$ be an arbitrary
non-negative initial data. Then, there exists $\gamma_{1}=\gamma
_{1}(\left\Vert \theta_{0}\right\Vert _{H^{\frac{3}{2}}})\in(1/2,1)$ such that
for each $\gamma\in\lbrack\gamma_{1},1)$ the IVP (\ref{1Dmodel}) has a unique
global (classical) $H^{\frac{3}{2}}$-solution.
\end{tma}

\begin{obs}
For each $\gamma\in(0,1)$, let $R_{\gamma}$ be the supremum of all $R>0$ such
that, for any$\;\theta_{0}\in H^{2}\;$with$\;\left\Vert \theta_{0}\right\Vert
_{\dot{H}^{\frac{3}{2}}(\mathbb{T})}^{1-\frac{2\gamma}{3}}\left\Vert
\theta_{0}\right\Vert _{L^{2}(\mathbb{T})}^{\frac{2\gamma}{3}}\leq R,$ the
unique $H^{\frac{3}{2}}$-solution of$\;$(\ref{1Dmodel}) with initial
data$\;\theta_{0}\;$does not blow up in finite time. In view of arguments in
the proof of Theorem \ref{global}, we have that $R_{\gamma}\rightarrow\infty$
as $\gamma\rightarrow1^{-}$. Moreover, in view of (\ref{rcondition}), we could
choose $\gamma_{1}$ depending on the quantity $\left\Vert \theta
_{0}\right\Vert _{\dot{H}^{\frac{3}{2}}(\mathbb{T})}^{1-\frac{2\gamma}{3}%
}\left\Vert \theta_{0}\right\Vert _{L^{2}(\mathbb{T})}^{\frac{2\gamma}{3}}$
(instead of $\left\Vert \theta_{0}\right\Vert _{H^{\frac{3}{2}}}$) which is
invariant by the scaling of (\ref{1Dmodel}).
\end{obs}

This paper is organized as follows. In Section 2 we recall some definitions,
notations and properties for Hilbert transform and fractional Laplacian
operator. Section 3 is devoted to the eventual regularity property. Finally,
in Section 4 we prove Theorem \ref{global}.

\section{Preliminaries}

Let $1\leq p\leq\infty$ and denote the norm of $L^{p}(\mathbb{T})$ by
$\left\Vert \cdot\right\Vert _{L^{p}}.$ For $s\in\mathbb{R}$ the norms of the
homogeneous Sobolev space $\dot{H}^{s}(\mathbb{T})$ and its nonhomogeneous
counterpart $H^{s}(\mathbb{T})$ are denoted by $\left\Vert \cdot\right\Vert
_{\dot{H}^{2}}=\left\Vert \Lambda^{s}\cdot\right\Vert _{L^{2}}$ and
$\left\Vert \cdot\right\Vert _{H^{2}}=\left\Vert \Lambda^{s}\cdot\right\Vert
_{L^{2}}+\left\Vert \cdot\right\Vert _{L^{2}}$, respectively. In turn, for
each $\alpha\in(0,1)$, the H\"{o}lder space $C^{\alpha}(\mathbb{T})$ is
endowed with the norm $\left\Vert \phi\right\Vert _{C^{\alpha}(\mathbb{T}%
)}=[\phi]_{C^{\alpha}(\mathbb{T})}+\left\Vert \phi\right\Vert _{L^{\infty
}(\mathbb{T})}$, where the seminorm $[\phi]$ is given by
\begin{equation}
\lbrack\phi]_{C^{\alpha}(\mathbb{T})}=\sup_{x,y\in\mathbb{T},\;x\neq y}%
\frac{|\phi(x)-\phi(y)|}{|x-y|^{\alpha}}. \label{seminorm1}%
\end{equation}

We recall that the periodic Hilbert transform $\mathcal{H}$ is defined by
means of Fourier transform as
\[
\widehat{\mathcal{H}\phi}(m)=-i\;\mbox{sign}(m)\widehat{\phi}(m)
\]
for all $m\in\mathbb{Z}_{\ast}$ and $\phi\in C^{\infty}(\mathbb{T}).$
Alternatively, in original variables we have the expression
\[
\mathcal{H}\phi(x)=\frac{1}{2\pi}P.V.\int_{\mathbb{T}}\frac{\phi(y)}%
{\tan(\frac{x-y}{2})}dy,
\]
which can be equivalently written as (see \cite{Calderon})
\begin{align}
\mathcal{H}\phi(x)  &  =\frac{1}{\pi}P.V.\int_{\mathbb{T}}\frac{\phi(y)}%
{x-y}dy+\frac{1}{\pi}\displaystyle\sum_{k\in\mathbb{Z}_{\ast}}\int
_{\mathbb{T}}\phi(y)\left(  \frac{1}{x-y-2k\pi}+\frac{1}{2k\pi}\right)
dy\nonumber\\
&  =\frac{1}{\pi}P.V.\int_{\mathbb{R}}\frac{\phi(y)}{x-y}dy. \label{hilbert2}%
\end{align}
In the last integral in (\ref{hilbert2}), recall that $P.V.$ is defined by
\begin{equation}
P.V.\int_{\mathbb{R}}\frac{\phi(y)}{x-y}dy=\lim_{\epsilon\rightarrow0}%
\int_{\epsilon<|x-y|<\frac{1}{\epsilon}}\frac{\phi(y)}{x-y}dy.\nonumber
\end{equation}
Hilbert transform commutates with derivatives and in particular we have that
\begin{equation}
\partial_{x}\mathcal{H}(\phi)(x)=\mathcal{H}(\partial_{x}\phi)(x).\nonumber
\end{equation}

For $0<\gamma<2$ and $\phi\in C^{\infty}(\mathbb{T}),$ the fractional
Laplacian $\Lambda^{\gamma}$ is defined by the following singular integral
(see \cite{Cordoba2} for more details)
\[
\Lambda^{\gamma}\phi(x)=c_{\gamma}\sum_{k\in\mathbb{Z}}\int_{\mathbb{T}}%
\frac{\phi(x)-\phi(x+y)}{|y-2\pi k|^{1+\gamma}}dy=c_{\gamma}P.V.\int
_{\mathbb{R}}\frac{\phi(x)-\phi(x+y)}{|y|^{1+\gamma}}dy
\]
where $c_{\gamma}$ is a normalization constant. For $\gamma\in(\gamma_{0},1],$
the constant $c_{\gamma}$ can be bounded from below and above by using
$\gamma_{0}$ and some universal constants $C.$ The exact expression of the
constant $c_{\gamma}$ is not necessary for our ends.

We recall the following commutator estimate for the fractional Laplacian (for
more details, see \cite{Ju} and references therein).

\begin{lem}
\label{comutadorkato} Suppose that $s>0$ and $1<p<\infty$. If $f,g\in
C^{\infty}(\mathbb{T})$ then
\[
\left\Vert \Lambda^{s}(fg)-f\Lambda^{s}g\right\Vert _{L^{p}(\mathbb{T})}\leq
C\left(  \left\Vert \partial_{x}f\right\Vert _{L^{p_{1}}(\mathbb{T}%
)}\left\Vert \Lambda^{s-1}g\right\Vert _{L^{p_{2}}(\mathbb{T})}+\left\Vert
\Lambda^{s}f\right\Vert _{L^{p_{3}}(\mathbb{T})}\left\Vert g\right\Vert
_{L^{p_{4}}(\mathbb{T})}\right)
\]
where $p_{1},p_{2},p_{3},p_{4}\in(1,\infty)$ are such that
\[
\frac{1}{p}=\frac{1}{p_{1}}+\frac{1}{p_{2}}=\frac{1}{p_{3}}+\frac{1}{p_{4}}%
\]
and $C>0$ is a constant depending on $s,p,p_{2}$ and $p_{3}$.
\end{lem}

The next lemma contains a property of the fractional Laplacian (see
\cite{Cordoba2} for more details and a proof in the two-dimensional case).

\begin{lem}
\label{laplacianofrac} Assume that $0<\gamma<2$ and $\phi\in C^{\infty
}(\mathbb{T})$. Then, we have the pointwise equality
\[
2\phi(x)\Lambda^{\gamma}\phi(x)=\Lambda^{\gamma}(\phi)^{2}(x)+D_{\gamma}%
(\phi)(x)\text{\ in }\mathbb{T}\text{, }%
\]
where
\[
D_{\gamma}(\phi)(x)=c_{\gamma}\sum_{k\in\mathbb{Z}}\int_{\mathbb{T}}%
\frac{(\phi(x)-\phi(x+y))^{2}}{|y-2\pi k|^{1+\gamma}}dy=c_{\gamma}%
P.V.\int_{\mathbb{R}}\frac{(\phi(x)-\phi(x+y))^{2}}{|y|^{1+\gamma}}dy.
\]

\end{lem}

We finish this section with a technical lemma that will be useful for our
purposes (see \cite[ Lemma B.1]{Constantin1}).

\begin{lem}
\label{lemaderivada}Let $\mathcal{K}\subset\mathbb{R}^{n}$ be compact and
$T>0$. Consider a function
\[
f:\mathcal{K}\times(0,T)\rightarrow\lbrack0,\infty)
\]
and assume that the functions
\[
f_{\lambda}(\cdot)=f(\lambda,\cdot):(0,T)\rightarrow\lbrack0,\infty
)\;\;\;\;\mbox{and}\;\;\;\;{f^{\prime}}_{\lambda}(\cdot)=(\partial
_{t}f)(\lambda,\cdot):(0,T)\rightarrow\mathbb{R}%
\]
are continuous, for each $\lambda\in\mathcal{K}.$ Additionally, assume that
the following properties hold true:

\begin{itemize}
\item[(i)] The families $\{f_{\lambda}\}_{\lambda\in\mathcal{K}}$ and
$\{f_{\lambda}^{\prime}\}_{\lambda\in\mathcal{K}}$ are uniformly
equicontinuous with respect to $t$;

\item[(ii)] For every $t\in(0,T)$, the functions
\[
f(\cdot,t):\mathcal{K}\rightarrow\lbrack0,\infty
)\;\;\;\;\mbox{and}\;\;\;\;(\partial_{t}f)(\cdot,t):\mathcal{K}\rightarrow
\mathbb{R}%
\]
are continuous. Moreover, define
\[
F(t)=\sup_{\lambda\in\mathcal{K}}f_{\lambda}(t).
\]
Then, for almost every $t\in(0,T)$ the function $F$ is differentiable at $t$
and there exists $\lambda_{\ast}=\lambda_{\ast}(t)\in\mathcal{K}$ such that
the following equalities hold simultaneously:
\[
F(t)=f_{\lambda_{\ast}}(t)\;\;\;\;\mbox{and}\;\;\;\;F^{\prime}(t)=f_{\lambda
_{\ast}}^{\prime}(t).
\]

\end{itemize}
\end{lem}

\section{Eventual regularization}

In this section, we show \textbf{a new} \textquotedblleft eventual
regularization\textquotedblright\;result for solutions of (\ref{1Dmodel}) in
which we provide an explicit control on the eventual regularity time $T^{\ast
}.$

Firstly, in \cite{Dong}, we can find the following local existence result for
(\ref{1Dmodel}): let $\gamma\in(0,1)$ and $\theta_{0}\in H^{\frac{3}{2}%
-\gamma}(\mathbb{T})$. Then, there exists $T>0$ such that the IVP
(\ref{1Dmodel}) has a unique strong solution
\[
\theta\in C([0,T);H^{\frac{3}{2}-\gamma}(\mathbb{T}))\cap L^{2}((0,T);H^{\frac
{3-\gamma}{2}}(\mathbb{T})).
\]
Using regularity techniques for the solution obtained in \cite{Dong} (see
Section 4 and estimates (\ref{local1})-(\ref{local7})), we can show that there
exists $0<T_{1}\leq T$ such that
\begin{equation}
\theta\in C([0,T_{1});H^{\frac{3}{2}}(\mathbb{T})) \label{Dong-local-2}%
\end{equation}
provided that $\theta_{0}\in H^{\frac{3}{2}}(\mathbb{T})$. Moreover,{%
\[
\theta\in C^{1}((0,T_{1});H^{\frac{1}{2}}(\mathbb{T})),\text{ }\partial
_{t}\theta\in L^{\infty}((0,\tilde{T});H^{\frac{1}{2}}(\mathbb{T}))\text{ and
}\theta\in L^{2}((0,\tilde{T});H^{\frac{3+\gamma}{2}}(\mathbb{T})),
\]
for all $0<\tilde{T}<T_{1}.$ }

In what follows, we state our \textquotedblleft eventual
regularization\textquotedblright\ result.

\begin{tma}
\label{eventual} Suppose that $\gamma\in\lbrack1/2,1)$ and $\theta_{0}\in
L^{\infty}(\mathbb{T})$ is non-negative. Let $\alpha\in(1-\gamma,1)$ and
define
\begin{equation}
T^{\ast}=C\alpha^{\frac{1}{1-\gamma}}\left\Vert \theta_{0}\right\Vert
_{L^{\infty}(\mathbb{T})}^{\frac{\gamma}{1-\gamma}}, \label{timeeventual}%
\end{equation}
where $C=\gamma^{-1}k_{1}k_{2}^{\frac{\gamma}{1-\gamma}}>0$ with $k_{1}$ and
$k_{2}$ being independent of $\alpha,\gamma$ and $\theta_{0}$. Let $\theta$ be
a solution of (\ref{1Dmodel}) in $C([0,T_{1});H^{\frac{3}{2}}(\mathbb{T}))$
with existence time $0<$ $T_{1}<\infty$. If $T^{\ast}<T_{1}$, then $\theta\in
C^{\infty}(\mathbb{T\times}(T^{\ast},T_{1}]).$
\end{tma}

\begin{obs}
Let us observe that the expression \textquotedblleft eventual
regularization\textquotedblright is used in the literature in the context of
weak solutions and $T_{1}=\infty.$ Nevertheless, in our range of $\gamma$, it
is not known whether (\ref{1Dmodel}) has global weak solution and then we need
to adapt this kind of result to our context but borrowing the same expression.
\end{obs}

In next lemma we recall a well-known result that assures that the control of
high-order Holder norms is sufficient to obtain smoothness. This result
essentially follows from \cite[Theorem 2.1]{Do} which extended the results of
\cite{Constantin2} about 2DQG to (\ref{1Dmodel}).

\begin{lem}
\label{lemmaholder}Let $\theta$ be a solution of (\ref{1Dmodel}) in the class
(\ref{Dong-local-2}) with non-negative initial data $\theta_{0}$. If
$0<t_{0}<t_{1}\leq T_{1}$ and
\begin{equation}
\theta\in L^{\infty}((t_{0},t_{1});C^{\alpha}(\mathbb{T}))\nonumber
\end{equation}
with $0<\gamma<1$ and $1-\gamma<\alpha$, then
\begin{equation}
\theta\in C^{\infty}(\mathbb{T\times}(t_{0},t_{1}]).\nonumber
\end{equation}

\end{lem}

\subsection{Proof of Theorem \ref{eventual}}

In view of Lemma \ref{lemmaholder}, we need to show that%
\[
\theta\in L^{\infty}((T^{\ast},T_{1});C^{\alpha}(\mathbb{T})),
\]
where $\alpha>1-\gamma$ and $T^{\ast}$ is as in (\ref{timeeventual}). For
that, we denote $\delta_{h}\theta(x,t)=\theta(x+h,t)-\theta(x,t)$ and define
\begin{equation}
v(x,t,h)=\frac{\delta_{h}\theta(x,t)}{(\xi^{2}(t)+|h|^{2})^{\frac{\alpha}{2}}%
}, \label{defv}%
\end{equation}
where $\xi:[0,\infty)\rightarrow\lbrack0,\infty)$ is a bounded decreasing
differentiable function which will be determined later. Notice that it is
sufficient to estimate $\left\Vert v(t)\right\Vert _{L^{\infty}}$ when
$\xi(t)=0$ in order to control the seminorm (\ref{seminorm1}) of $\theta$ in
$C^{\alpha}(\mathbb{T}).$

We start by providing estimates for $Lv^{2}$ where $L$ is the operator of the
corresponding equation satisfied by $\delta_{h}\theta$. We split the proof
into a sequence of lemmas. Taking the differences in (\ref{1Dmodel}) evaluated
in $x+h$ and $x,$ it follows that
\begin{equation}
\partial_{t}\delta_{h}\theta-\mathcal{H}\theta\partial_{x}\delta_{h}%
\theta-\delta_{h}\mathcal{H}\theta\partial_{h}\delta_{h}\theta+\Lambda
^{\gamma}\delta_{h}\theta=0, \label{L1}%
\end{equation}
which gives $L=\partial_{t}-u\partial_{x}-\delta_{h}u\partial_{h}%
+\Lambda^{\gamma}$with $u=\mathcal{H}\theta$. Combining (\ref{L1}) and Lemma
\ref{laplacianofrac} we obtain that $v^{2}$ satisfies
\begin{equation}
Lv^{2}+\frac{1}{(\xi^{2}(t)+|h|^{2})^{\alpha}}D_{\gamma}(\delta_{h}%
\theta)=-2\alpha\xi^{\prime}\frac{\xi}{\xi^{2}+|h|^{2}}v^{2}+2\alpha\frac
{h}{\xi^{2}+|h|^{2}}\delta_{h}\mathcal{H}\theta v^{2}. \label{L2}%
\end{equation}

Next we estimate the term $D_{\gamma}(\delta_{h}\theta)$.

\begin{lem}
\label{lemanlinear} Let $0<\gamma_{0}\leq\gamma<1$ and $\alpha\in(1-\gamma
,1)$. Then, there exists a positive constant $c_{0}=c_{0}(\gamma_{0})$ such
that
\[
\displaystyle D_{\gamma}(\delta_{h}\theta)(x)\geq\frac{1}{c_{0}|h|^{\gamma}%
}\left(  \frac{|v(x,h)|}{\left\Vert v\right\Vert _{L^{\infty}}}\right)
^{\frac{\gamma}{1-\alpha}}|\delta_{h}\theta(x)|^{2},
\]
for all $x,h\in\mathbb{T}.$
\end{lem}

\textbf{Proof.} Let $\chi$ be a smooth radially non-decreasing cutoff function
such that $\chi$ vanishes for $|x|\leq1$, $\chi$ is identically $1$ for
$|x|\geq2$ and its derivative verifies $|\chi^{\prime}|\leq2$. For $R\geq
4|h|$, we can estimate
\begin{align}
D_{\gamma}(\delta_{h}\theta)(x)  &  \geq c_{\gamma}\int_{\mathbb{R}}%
\frac{(\delta_{h}\theta(x)-\delta_{h}\theta(x+y))^{2}}{|y|^{1+\gamma}}%
\chi\left(  \frac{|y|}{R}\right)  dy\nonumber\\
&  \geq c_{\gamma}|\delta_{h}\theta(x)|^{2}\int_{|y|\geq2R}\frac
{1}{|y|^{1+\gamma}}dy-2c_{\gamma}|\delta_{h}\theta(x)|\left\vert
\int_{\mathbb{R}}\frac{\delta_{h}\theta(x+y)}{|y|^{1+\gamma}}\chi\left(
\frac{|y|}{R}\right)  dy\right\vert \nonumber\\
&  \geq c_{\gamma}\frac{|\delta_{h}\theta(x)|^{2}}{R^{\gamma}}-2c_{\gamma
}|\delta_{h}\theta(x)|\left\vert \int_{\mathbb{R}}(\theta(x+y)-\theta
(x))\delta_{-h}\left(  \frac{\chi(\frac{|y|}{R})}{|y|^{1+\gamma}}\right)
dy\right\vert . \label{lemma31estimativa11}%
\end{align}

Denoting $g(y)=\chi\left(  \frac{|y|}{R}\right)  |y|^{-(1+\gamma)}$ and using
mean value theorem, we obtain
\[
|\delta_{-h}g(y)|\leq|h|\max_{0\leq\lambda\leq1}|g^{\prime}(y-\lambda h)|\leq
c_{1}\frac{|h|}{|y|^{2+\gamma}}\mathbf{1}_{\left\{  \frac{3R}{4}%
\leq|y|\right\}  },
\]
for some constant $c_{1}\geq1$. Hence the integral in
(\ref{lemma31estimativa11}) can be estimated as
\begin{align}
\left\vert \int_{\mathbb{R}}(\theta(x+y)-\theta(x))\delta_{-h}\left(
\frac{\chi(\frac{|y|}{R})}{|y|^{1+\gamma}}\right)  dy\right\vert  &  \leq
c_{1}|h|\int_{|y|\geq\frac{3R}{4}}\frac{|\delta_{y}\theta(x)|}{(\xi
^{2}+|y|^{2})^{\frac{\alpha}{2}}}\frac{(\xi^{2}+|y|^{2})^{\frac{\alpha}{2}}%
}{|y|^{2+\gamma}}dy\nonumber\\
&  \leq c_{1}|h|\left\Vert v\right\Vert _{L^{\infty}}\int_{|y|\geq\frac{3R}%
{4}}\frac{(\xi^{2}+|y|^{2})^{\frac{\alpha}{2}}}{|y|^{2+\gamma}}dy\nonumber\\
&  \leq c_{1}c_{2}\frac{|h|\left\Vert v\right\Vert _{L^{\infty}}}{R^{\gamma}%
}\left(  \frac{\xi^{\alpha}}{R}+\frac{1}{R^{1-\alpha}}\right)  ,
\label{lemma31estimativa2}%
\end{align}
for some constant $c_{2}\geq1$.

Now choose $R>0$ defined by
\begin{equation}
R=\left[  \frac{8c_{1}c_{2}\left\Vert v\right\Vert _{L^{\infty}}}%
{|v(x,t)|}\right]  ^{\frac{1}{1-\alpha}}|h|. \label{lemma31defR}%
\end{equation}
Since $c_{1}c_{2}\geq1$ and $|v(x,t)|\leq\left\Vert v\right\Vert _{L^{\infty}%
},$ it follows that $R\geq8^{\frac{1}{1-\alpha}}|h|\geq4|h|$.

Using (\ref{lemma31defR}), we can rewrite estimate (\ref{lemma31estimativa2})
as
\begin{align}
\left\vert \int_{\mathbb{R}}(\theta(x+y)-\theta(x))\delta_{-h}\left(
\frac{\chi(\frac{|y|}{R})}{|y|^{1+\gamma}}\right)  dy\right\vert  &  \leq
\frac{1}{R^{\gamma}}\left[  \frac{c_{1}c_{2}\xi^{\alpha}\left\Vert
v\right\Vert _{L^{\infty}}}{(8c_{1}c_{2})^{\frac{1}{1-\alpha}}}\left[
\frac{|v(x,t)|}{\left\Vert v\right\Vert _{L^{\infty}}}\right]  ^{\frac
{1}{1-\alpha}}+\frac{|v(x,t)||h|^{\alpha}}{8}\right] \nonumber\\
&  \leq\frac{|\xi|^{\alpha}+|h|^{\alpha}}{8R^{\gamma}}v(x,t)|\nonumber\\
&  \leq\frac{|\delta_{h}\theta(x)|}{4R^{\gamma}}, \label{lemma31estimativa3}%
\end{align}
because $a^{\frac{1}{1-\gamma}}\leq a$ for all $a\leq1$.

Combining estimates (\ref{lemma31estimativa11}) and (\ref{lemma31estimativa3})
and using the definition of $R,$ we arrive at
\begin{align}
D_{\gamma}(\delta_{h}\theta)(x)  &  \geq c_{\gamma}\frac{|\delta_{h}%
\theta(x)|^{2}}{2R^{\gamma}}\nonumber\\
&  \geq c_{\gamma}\frac{1}{2(8c_{1}c_{2})^{\frac{\gamma}{1-\alpha}}%
|h|^{\gamma}}\left[  \frac{|v(x,t)|}{\left\Vert v\right\Vert _{L^{\infty}}%
}\right]  ^{\frac{\gamma}{1-\alpha}}|\delta_{h}\theta(x)|^{2},
\label{lema31estimativafinal}%
\end{align}
as required.\fin

\begin{obs}
Notice that the condition $0<\gamma_{0}\leq\gamma<1$ arises from the need of
controlling terms with the power $\frac{\gamma}{1-\alpha}$ in
(\ref{lemma31defR}) and (\ref{lema31estimativafinal}).
\end{obs}

In the next lemma, we define $\xi$ by an ordinary differential equation and
obtain an estimate for the first term on the right hand side of (\ref{L2}).

\begin{lem}
\label{lemaH1}Let $\gamma_{0}>0,$ $\gamma\in\lbrack\gamma_{0},1)$ and
$\alpha\in(1-\gamma,1)$. There exists a positive constant $k_{1}=k_{1}%
(\gamma_{0})$ such that if
\begin{equation}
\xi^{\prime}=-\frac{1}{\alpha k_{1}}\xi^{1-\gamma}, \label{lemma311}%
\end{equation}
then
\begin{equation}
-2\alpha\xi^{\prime}\frac{\xi}{\xi^{2}+|h|^{2}}v^{2}\leq\frac{1}%
{8c_{0}|h|^{\gamma}}v^{2},\text{ for all }x,h\in\mathbb{T}, \label{lemma312}%
\end{equation}
where $c_{0}$ is the same constant appearing in Lemma \ref{lemanlinear}.
\end{lem}

\textbf{Proof.} Substituting (\ref{lemma311}) on the left hand side of
(\ref{lemma312}), we conclude that
\begin{align*}
-2\alpha\xi^{\prime}\frac{\xi}{\xi^{2}+|h|^{2}}v^{2}  &  \leq\frac
{2(k_{1})^{-1}\xi^{2-\gamma}}{\xi^{2}+|h|^{2}}v^{2}\\
&  \leq\frac{2(k_{1})^{-1}(\xi^{2}+|h|^{2})^{1-\frac{\gamma}{2}}}{\xi
^{2}+|h|^{2}}v^{2}\\
&  \leq\frac{2(k_{1})^{-1}}{(\xi^{2}+|h|^{2})^{\frac{\gamma}{2}}}v^{2}%
\leq\frac{2(k_{1})^{-1}}{|h|^{\gamma}}v^{2}.
\end{align*}
Just choosing $k_{1}=16c_{0}$, we are done. \fin

Now we need to estimate the term in (\ref{L2}) that depends on the Hilbert
transform. We do that in the next two lemmas. First, we work with the factor
$\delta_{h}\mathcal{H}\theta(x)$ in the second term of the right hand side of
(\ref{L2}).

\begin{lem}
\label{lemaH2}Let $\gamma\in[\gamma_{0},1)$ and $\alpha\in(1-\gamma,1)$ . If
$\rho\geq4|h|$, then
\begin{equation}
|\delta_{h}\mathcal{H}\theta(x)|\leq C\left[  \rho^{\frac{\gamma}{2}%
}(D_{\gamma}(\delta_{h}\theta)(x))^{\frac{1}{2}}+\left\Vert v\right\Vert
_{L^{\infty}}\left(  \frac{|h|\xi^{\alpha}}{\rho}+\frac{|h|}{\rho^{1-\alpha}%
}\right)  \right]  , \label{Hestimative}%
\end{equation}
for all $x,h\in\mathbb{T}$.
\end{lem}

\textbf{Proof.} Let $\chi$ be a smooth radially non-decreasing cutoff function
that vanishes for $|x|\leq1$ and is equal $1$ for $|x|\geq2,$ and such that
the derivative satisfies $|\chi^{\prime}|\leq2$. From (\ref{hilbert2}), we
obtain
\[
\delta_{h}(\mathcal{H}\theta)(x)=\mathcal{H}(\delta_{h}\theta)(x)=-\frac
{1}{\pi}P.V\int_{\mathbb{R}}\frac{\delta_{h}\theta(x+y)}{y}dy.
\]
For $\epsilon>0,$ it follows that
\begin{align}
\int_{\epsilon\leq|y|\leq\frac{1}{\epsilon}}\frac{\delta_{h}\theta(x+y)}{y}dy
&  =\int_{\epsilon\leq|y|\leq\frac{1}{\epsilon}}\frac{\delta_{h}%
\theta(x+y)-\delta_{h}\theta(x)}{y}dy\nonumber\\
&  :=I+J, \label{Hestimative1}%
\end{align}
where
\[
I=\int_{\epsilon\leq|y|\leq\frac{1}{\epsilon}}\left[  1-\chi\left(  \frac
{|y|}{\rho}\right)  \right]  \frac{\delta_{h}\theta(x+y)-\delta_{h}\theta
(x)}{y}dy
\]
and
\[
J=\int_{\epsilon\leq|y|\leq\frac{1}{\epsilon}}\chi\left(  \frac{|y|}{\rho
}\right)  \frac{\delta_{h}\theta(x+y)-\delta_{h}\theta(x)}{y}dy.
\]

Applying H\"{o}lder inequality and taking $\epsilon\leq\frac{1}{2\rho}$, we
can estimate $I$ as follows%
\begin{align}
|I|  &  \leq\int_{\epsilon\leq|y|\leq2\rho}\frac{|\delta_{h}\theta
(x+y)-\delta_{h}\theta(x)|}{|y|}dy\nonumber\\
&  \leq\left(  \int_{|y|\leq2\rho}\frac{1}{|y|^{1-\gamma}}dy\right)
^{\frac{1}{2}}\left(  \int_{\epsilon\leq|y|\leq2\rho}\frac{|\delta_{h}%
\theta(x+y)-\delta_{h}\theta(x)|^{2}}{|y|^{1+\gamma}}dy\right)  ^{\frac{1}{2}%
}\nonumber\\
&  \leq C\rho^{\frac{\gamma}{2}}\left(  \int_{\epsilon\leq|y|}\frac
{|\delta_{h}\theta(x+y)-\delta_{h}\theta(x)|^{2}}{|y|^{1+\gamma}}dy\right)
^{\frac{1}{2}}. \label{HestimativeI}%
\end{align}

For $J$, we have that
\begin{equation}
|J|\leq\int_{\epsilon\leq|y|\leq\frac{1}{\epsilon}}\left\vert \delta
_{-h}\left(  \frac{\chi\left(  \frac{|y|}{\rho}\right)  }{y}\right)
\right\vert \text{ }\left\vert (\theta(x+y)-\theta(x))\right\vert dy.
\label{estimativaJ1}%
\end{equation}
Taking $g(y)=\chi\left(  \frac{|y|}{\rho}\right)  y^{-1}$ and applying mean
value theorem, we arrive at
\begin{equation}
|\delta_{-h}g(y)|\leq|h|\max_{0\leq\lambda\leq1}|g^{\prime}(y-\lambda h)|\leq
C|h|\frac{\mathbf{1}_{\left\{  \frac{3\rho}{4}\leq|y|\right\}  }}{|y|^{2}}.
\label{estimativaJ2}%
\end{equation}
Substituting (\ref{estimativaJ2}) into (\ref{estimativaJ1}) and taking
$\epsilon\leq\frac{3\rho}{4}$, we obtain
\begin{align}
J  &  \leq C|h|\int_{\frac{3\rho}{4}\leq|y|\leq\frac{1}{\epsilon}}%
\frac{|(\theta(x+y)-\theta(x))|}{|y|^{2}}dy\nonumber\\
&  \leq C|h|\left\Vert v\right\Vert _{L^{\infty}}\int_{\frac{3\rho}{4}\leq
|y|}\frac{(\xi^{2}+|y|^{2})^{\frac{\alpha}{2}}}{|y|^{2}}dy\nonumber\\
&  \leq C|h|\left\Vert v\right\Vert _{L^{\infty}}\left(  \frac{\xi^{\alpha}%
}{\rho}+\frac{1}{\rho^{1-\alpha}}\right)  . \label{HestimativeJ}%
\end{align}

The estimate (\ref{Hestimative}) follows by inserting (\ref{HestimativeI}) and
(\ref{HestimativeJ}) in (\ref{Hestimative1}) and making $\epsilon\rightarrow
0$. \fin

Next we provide a condition on the initial data of $\xi$ to relate
(\ref{Hestimative}) to the estimates of the other terms.

\begin{lem}
\label{lemaH3} Let $\gamma\in[\gamma_{0},1)$ and $\alpha\in(1-\gamma,1),$ and
assume that
\begin{equation}
\left\Vert v\right\Vert _{L^{\infty}}\leq\frac{4\left\Vert \theta
_{0}\right\Vert _{L^{\infty}}}{\xi_{0}^{\alpha}}. \label{lemma4condicaoM}%
\end{equation}
There exists a constant $k_{2}=k_{2}(\gamma_{0})\geq1$ such that if
\[
\xi_{0}=(k_{2}\alpha\left\Vert \theta_{0}\right\Vert _{L^{\infty}})^{\frac
{1}{1-\gamma}},
\]
then
\[
2\alpha\frac{|h|}{\xi^{2}+|h|^{2}}|\delta_{h}\mathcal{H}\theta|v^{2}\leq
\frac{1}{2(\xi^{2}+|h|^{2})^{\alpha}}D_{\gamma}(\delta_{h}\theta)(x)+\frac
{1}{8c_{0}|h|^{\gamma}}v^{2},
\]
for all $x,h\in\mathbb{T}$ with $|h|\leq\xi_{0}.$
\end{lem}

\textbf{Proof.} From (\ref{Hestimative}) and Young's inequality for products
it follows that
\begin{align}
2\alpha\frac{|h|}{\xi^{2}+|h|^{2}}|\delta_{h}\mathcal{H}\theta|v^{2}  &
\leq\frac{|h|}{2(\xi^{2}+|h|^{2})^{\alpha}}D_{\gamma}(\delta_{h}%
\theta)(x)\nonumber\\
&  +C\alpha\frac{|h|^{2}}{\xi^{2}+|h|^{2}}\left[  \frac{\alpha\rho^{\gamma
}v^{2}}{(\xi^{2}+|h|^{2})^{1-\alpha}}+\left\Vert v\right\Vert _{L^{\infty}%
}\left(  \frac{\xi^{\alpha}}{\rho}+\frac{1}{\rho^{1-\alpha}}\right)  \right]
v^{2}.\nonumber
\end{align}

It is sufficient to show that
\begin{equation}
C\alpha\left[  \frac{\alpha\rho^{\gamma}v^{2}}{(\xi^{2}+|h|^{2})^{1-\alpha}%
}+\left\Vert v\right\Vert _{L^{\infty}}\left(  \frac{\xi^{\alpha}}{\rho}%
+\frac{1}{\rho^{1-\alpha}}\right)  \right]  \leq\frac{1}{8c_{0}|h|^{\gamma}}.
\label{lemma41}%
\end{equation}

Choose
\begin{equation}
\rho=4(\xi^{2}+|h|^{2})^{\frac{1}{2}}, \label{lema4rho}%
\end{equation}
and note that $\rho\geq4|h|$. Combining (\ref{lemma4condicaoM}) and
(\ref{lema4rho}), we obtain
\begin{equation}
\frac{\alpha\rho^{\gamma}v^{2}}{(\xi^{2}+|h|^{2})^{1-\alpha}}\leq
C\frac{\alpha\left\Vert \theta_{0}\right\Vert _{L^{\infty}}^{2}(\xi
^{2}+|h|^{2})^{\frac{\gamma}{2}}}{\xi_{0}^{2\alpha}(\xi^{2}+|h|^{2}%
)^{1-\alpha}}. \label{lema411}%
\end{equation}
For $|h|\leq\xi_{0}$, $1-\gamma<\alpha<1$ and $\xi\leq\xi_{0}$, we deduce
that
\begin{equation}
\frac{(\xi^{2}+|h|^{2})^{\frac{\gamma}{2}}}{\xi_{0}^{2\alpha}(\xi^{2}%
+|h|^{2})^{1-\alpha}}\leq C\frac{(\xi^{2}+|h|^{2})^{\gamma+\alpha-1}}{\xi
_{0}^{2\alpha}(\xi^{2}+|h|^{2})^{\frac{\gamma}{2}}}\leq C\frac{1}{\xi
_{0}^{2(1-\gamma)}|h|^{\gamma}}. \label{lema412}%
\end{equation}
Estimates (\ref{lema411}) and (\ref{lema412}) yield
\begin{equation}
\frac{\alpha\rho^{\gamma}v^{2}}{(\xi^{2}+|h|^{2})^{1-\alpha}}\leq
C\frac{\alpha\left\Vert \theta_{0}\right\Vert _{L^{\infty}}^{2}}{\xi
_{0}^{2(1-\gamma)}|h|^{\gamma}}. \label{lema413}%
\end{equation}

Proceeding similarly, one also can show that
\begin{align}
\frac{\left\Vert v\right\Vert _{L^{\infty}}\xi^{\alpha}}{\rho}+\frac
{\left\Vert v\right\Vert _{L^{\infty}}}{\rho^{1-\alpha}}  &  \leq
C\frac{\left\Vert \theta_{0}\right\Vert _{L^{\infty}}}{\xi_{0}^{\alpha}%
}\left(  \frac{\xi^{\alpha}}{(\xi^{2}+|h|^{2})^{\frac{1}{2}}}+\frac{1}%
{(\xi^{2}+|h|^{2})^{\frac{1-\alpha}{2}}}\right) \nonumber\\
&  \leq C\frac{\left\Vert \theta_{0}\right\Vert _{L^{\infty}}}{\xi_{0}%
^{\alpha}}\frac{(\xi^{2}+|h|^{2})^{\frac{\gamma+\alpha-1}{2}}}{|h|^{\gamma}%
}\nonumber\\
&  \leq C\frac{\left\Vert \theta_{0}\right\Vert _{L^{\infty}}}{\xi
_{0}^{1-\gamma}|h|^{\gamma}}. \label{lema414}%
\end{align}

Adding (\ref{lema413}) and (\ref{lema414}), we get
\begin{equation}
\frac{\alpha\rho^{\gamma}v^{2}}{(\xi^{2}+|h|^{2})^{1-\alpha}}+\left\Vert
v\right\Vert _{L^{\infty}}\left(  \frac{\xi^{\alpha}}{\rho}+\frac{1}%
{\rho^{1-\alpha}}\right)  \leq C\left(  \frac{\alpha\left\Vert \theta
_{0}\right\Vert _{L^{\infty}}^{2}}{\xi_{0}^{2(1-\gamma)}}+\frac{\left\Vert
\theta_{0}\right\Vert _{L^{\infty}}}{\xi_{0}^{1-\gamma}}\right)  \frac
{1}{|h|^{\gamma}}. \label{lemma4est1}%
\end{equation}
In view of (\ref{lemma4est1}), notice that we only need to show that
\begin{equation}
C\alpha\left(  \frac{\alpha\left\Vert \theta_{0}\right\Vert _{L^{\infty}}^{2}%
}{\xi_{0}^{2(1-\gamma)}}+\frac{\left\Vert \theta_{0}\right\Vert _{L^{\infty}}%
}{\xi_{0}^{1-\gamma}}\right)  \leq\frac{1}{8c_{0}}. \label{lemma4est2}%
\end{equation}
For that, we simply choose $\xi_{0}$ satisfying
\[
\frac{\left\Vert \theta_{0}\right\Vert _{L^{\infty}}}{\xi_{0}^{1-\gamma}}%
\leq\frac{1}{16\tilde{C}c_{0}\alpha},
\]
for some constant $\tilde{C}$ such that $\tilde{C}\geq C$ and $\tilde{C}%
c_{0}\geq1$. Combining (\ref{lemma4est1}) and (\ref{lemma4est2}), we conclude
(\ref{lemma41}). \fin

Finally, we are ready to conclude the proof of Theorem \ref{eventual}. \ For
$1-\gamma<\alpha\leq1/2$, let
\[
M=\frac{4\left\Vert \theta_{0}\right\Vert _{L^{\infty}}}{\xi_{0}^{\alpha}}%
\]
and define
\begin{equation}
t_{\ast}=\sup\{0\leq t<T_{1}:\left\Vert v(\tau)\right\Vert _{L^{\infty}%
}<M\text{ for all }\tau\in\lbrack0,t]\}. \label{testrela}%
\end{equation}
Note that $t_{\ast}$ is well-defined since (\ref{defv}) provides
\[
\left\Vert v(0)\right\Vert _{L^{\infty}}\leq\frac{M}{2}.
\]

We are going to show that $t_{\ast}=T_{1}$. Suppose by contradiction that
$t_{\ast}<T_{1}$. Since $v$ is continuous and periodic in $x$ and $h$ there
exists $(x_{0},h_{0})\in\mathbb{T}\times\mathbb{T}$ such that $|v(x_{0}%
,t_{\ast},h_{0})|=M$. We claim that $|h_{0}|\leq\xi_{0}$. Indeed, if
$|h_{0}|\geq\xi_{0}$, then
\[
|v(x_{0},t_{\ast},h_{0})|\leq\frac{2\left\Vert \theta\right\Vert _{L^{\infty}%
}}{|h_{0}|^{\alpha}}\leq\frac{2\left\Vert \theta_{0}\right\Vert _{L^{\infty}}%
}{\xi_{0}^{\alpha}}\leq\frac{M}{2}.
\]

Applying Lemmas \ref{lemaH1} and \ref{lemaH3} in (\ref{L2}) for $t\in
(0,t_{\ast}]$, we obtain the estimate
\begin{equation}
Lv^{2}+\frac{1}{(\xi^{2}+|h|^{2})^{\alpha}}D_{\gamma}(\delta_{h}\theta
)\leq\frac{1}{2(\xi^{2}+|h|^{2})^{\alpha}}D_{\gamma}(\delta_{h}\theta
)+\frac{1}{4c_{0}|h|^{\gamma}}v^{2}, \label{L4}%
\end{equation}
for all $x,h\in\mathbb{T}$ with $|h|\leq\xi_{0}$. Using Lemma
\ref{lemanlinear} we can rewrite (\ref{L4}) as follows
\begin{equation}
Lv^{2}+\frac{1}{4c_{0}|h|^{\gamma}}\left[  \left(  \frac{|v(x,h)|}{\left\Vert
v\right\Vert _{L^{\infty}}}\right)  ^{\frac{\gamma}{1-\alpha}}-1\right]
v^{2}+\frac{1}{4c_{0}|h|^{\gamma}}\left(  \frac{|v(x,h)|}{\left\Vert
v\right\Vert _{L^{\infty}}}\right)  ^{\frac{\gamma}{1-\alpha}}v^{2}\leq0.
\label{L6}%
\end{equation}

Now consider $\epsilon>0$ such that
\begin{equation}
\left\Vert v(t)\right\Vert _{L^{\infty}}\geq\frac{7M}{8},\text{ for all }%
t\in\lbrack t_{\ast}-\epsilon,t_{\ast}). \label{esti7M8}%
\end{equation}
Given $t\in\lbrack t_{\ast}-\epsilon,t_{\ast})$, consider $(x_{t},h_{t}%
)\in\mathbb{T}\times\mathbb{T}$ such that the function $(x,h)\rightarrow
v^{2}(x,t,h)$ reaches its maximum. At this point, we have that $\partial
_{x}v^{2}=\partial_{h}v^{2}=0$, $\Lambda^{\gamma}v^{2}\geq0$ and $|h_{t}%
|\leq\xi_{0}$, which leads us to
\begin{equation}
(\partial_{t}v^{2})(x_{t},t,h_{t})\leq Lv^{2}(x_{t},t,h_{t}).
\label{partialt1}%
\end{equation}
Using (\ref{esti7M8}) and $|h_{t}|\leq\xi_{0}$, we deduce that
\begin{equation}
\frac{49M}{256c_{0}\xi_{0}^{\gamma}}\leq\frac{v^{2}(x_{t},t,h_{t})}%
{4c_{0}|h_{t}|^{\gamma}}. \label{partialt2}%
\end{equation}
Next, adding (\ref{partialt1}) and (\ref{partialt2}), we conclude
\begin{equation}
(\partial_{t}v^{2})(x_{t},t,h_{t})+\frac{49M}{256c_{0}\xi_{0}^{\gamma}}\leq
Lv^{2}(x_{t},t,h_{t})+\frac{v^{2}(x_{t},t,h_{t})}{4c_{0}|h_{t}|^{\gamma}}.
\label{partialt3}%
\end{equation}
Combining estimate (\ref{L6}) at the point $(x_{t},t,h_{t})$ with
(\ref{partialt3}) and using that $|v(x_{t},t,h_{t})|=\left\Vert
v(t)\right\Vert _{L^{\infty}},$ it follows that
\[
(\partial_{t}v^{2})(x_{t},t,h_{t})\leq-\frac{49M}{256c_{0}\xi_{0}^{\gamma}},
\]
for all $t\in\lbrack t_{\ast}-\epsilon,t_{\ast})$.

Lemma \ref{lemaderivada} with $f(t,\lambda)=v(x,t,h)^{2}$ and $\lambda
=(x,h)\in\mathcal{K}=\mathbb{T}\times\mathbb{T}$ yields
\begin{equation}
\frac{d}{dt}\left\Vert v(t)\right\Vert _{L^{\infty}}^{2}\leq(\partial_{t}%
v^{2})(x_{t},t,h_{t})\leq-\frac{49M}{256c_{0}\xi_{0}^{\gamma}},
\label{aux-est-1000}%
\end{equation}
for all $t\in\lbrack t_{\ast}-\epsilon,t_{\ast})$. Integrating
(\ref{aux-est-1000}), we arrive at
\[
\left\Vert v(t_{\ast})\right\Vert _{L^{\infty}}<M,
\]
which contradicts (\ref{testrela}). Consequently $t_{\ast}=T_{1}$ and $v\in
L^{\infty}(\mathbb{T\times}(0,T_{1})).$

Next, taking $\xi_{0}=(k_{2}\alpha\left\Vert \theta_{0}\right\Vert
_{L^{\infty}})^{\frac{1}{1-\gamma}}$ as in Lemma \ref{lemaH3}, then the
solution of (\ref{lemma311}) is given by
\begin{equation}
\xi(t)=\left\{
\begin{array}
[c]{rc}%
\lbrack\xi_{0}^{\gamma}-\frac{\gamma}{\alpha k_{1}}t]^{\frac{1}{\gamma}}, &
\mbox{if}\text{ \ }0\leq t\leq T_{\ast}\\
0, & \mbox{if}\text{ \ }T_{\ast}<t<T_{1}%
\end{array}
\right.  ,\nonumber
\end{equation}
where
\[
T_{\ast}=C\alpha^{\frac{1}{1-\gamma}}\left\Vert \theta_{0}\right\Vert
_{L^{\infty}}^{\frac{\gamma}{1-\gamma}}\text{ with }C=\gamma^{-1}k_{1}%
k_{2}^{\frac{\gamma}{1-\gamma}}.
\]
Since $\xi(t)=0$ for $T_{\ast}<t<T_{1}$, it follows that
\[
\lbrack\theta(\cdot,t)]_{C^{\alpha}}\leq C\left\Vert v(t)\right\Vert
_{L^{\infty}}\leq M,\text{ for all }T_{\ast}<t<T_{1},
\]
and we are done.\fin

\section{Proof of Theorem \ref{global}}

Firstly we obtain an explicit lower bound of the local existence time with
$H^{\frac{3}{2}}$-initial data. For that, we need an \textit{a priori}
estimate of $H^{\frac{3}{2}}$-norm of the solution and after compare with the
eventual regularization time $T^{\ast}$ (\ref{timeeventual}).

Formally applying $\Lambda^{\frac{3}{2}}$ in (\ref{1Dmodel}) and then multiply
by $\Lambda^{\frac{3}{2}}\theta,$ we obtain
\begin{align}
\frac{1}{2}\frac{d}{dt}\left\Vert \Lambda^{\frac{3}{2}}\theta(t)\right\Vert
_{L^{2}}^{2}+\left\Vert \Lambda^{\frac{3+\gamma}{2}}\theta(t)\right\Vert
_{L^{2}}^{2}=  &  \int\Lambda^{\frac{3}{2}}\theta\Lambda^{\frac{3}{2}%
}\big(\mathcal{H}\theta\theta_{x}\big)(t,x)dx\nonumber\\
=  &  I_{1}+I_{2}, \label{local1}%
\end{align}
where
\[
I_{1}=\int\Lambda^{\frac{3}{2}}\theta\Big[\Lambda^{\frac{3}{2}}%
\big(\mathcal{H}\theta\theta_{x}\big)-\mathcal{H}\theta\Lambda^{\frac{3}{2}%
}\theta_{x}\Big](t,x)dx
\]
and
\[
I_{2}=\int\Lambda^{\frac{3}{2}}\theta\mathcal{H}\theta\Lambda^{\frac{3}{2}%
}\theta_{x}(t,x)dx.\;\;\;\;\;\;\;\;\;\;\;\;\;\;\;\;\;\;\;\;\;\;\;
\]
Using H\"{o}lder inequality and Lemma \ref{comutadorkato} with $p=2$,
$p_{1}=p_{4}=\frac{3}{\gamma}$ and $p_{2}=p_{3}=\frac{6}{3-2\gamma}$, we
conclude
\begin{align}
I_{1}  &  \leq\left\Vert \Lambda^{\frac{3}{2}}\theta(t)\right\Vert _{L^{2}%
}\left\Vert \Lambda^{\frac{3}{2}}\big(\mathcal{H}\theta\theta_{x}%
\big)(t)-\mathcal{H}\theta\Lambda^{\frac{3}{2}}\theta_{x}(t)\right\Vert
_{L^{2}}\nonumber\\
&  \leq C\left\Vert \Lambda^{\frac{3}{2}}\theta(t)\right\Vert _{L^{2}%
}\left\Vert \Lambda\theta(t)\right\Vert _{L^{\frac{3}{\gamma}}}\left\Vert
\Lambda^{\frac{3}{2}}\theta(t)\right\Vert _{L^{\frac{6}{3-2\gamma}}}.
\label{local2}%
\end{align}
An integration by parts and H\"{o}lder inequality lead us to
\[
I_{2}\leq\left\Vert \Lambda^{\frac{3}{2}}\theta(t)\right\Vert _{L^{2}%
}\left\Vert \Lambda\theta(t)\right\Vert _{L^{\frac{3}{\gamma}}}\left\Vert
\Lambda^{\frac{3}{2}}\theta(t)\right\Vert _{L^{\frac{6}{3-2\gamma}}}.
\]
Applying Gagliardo-Nirenberg inequality and (\ref{max-L2}), we conclude that
\[
\left\Vert \Lambda\theta(t)\right\Vert _{L^{\frac{3}{\gamma}}}\leq C\left\Vert
\Lambda^{\frac{3}{2}}\theta(t)\right\Vert _{L^{2}}^{1-\frac{2\gamma}{9}%
}\left\Vert \theta_{0}\right\Vert _{L^{2}}^{\frac{2\gamma}{9}}%
\]
and
\begin{equation}
\left\Vert \Lambda^{\frac{3}{2}}\theta(t)\right\Vert _{L^{\frac{6}{3-2\gamma}%
}}\leq C\left\Vert \Lambda^{\frac{3+\gamma}{2}}\theta(t)\right\Vert _{L^{2}%
}^{\frac{9+2\gamma}{3(3+\gamma)}}\left\Vert \theta_{0}\right\Vert _{L^{2}%
}^{\frac{\gamma}{3(3+\gamma)}}. \label{local5}%
\end{equation}
Thus, using (\ref{local2})-(\ref{local5}), we can estimate the right hand side
of (\ref{local1}) as
\begin{align}
\frac{1}{2}\frac{d}{dt}\left\Vert \Lambda^{\frac{3}{2}}\theta(t)\right\Vert
_{L^{2}}^{2}+  &  \left\Vert \Lambda^{\frac{3+\gamma}{2}}\theta(t)\right\Vert
_{L^{2}}^{2}\nonumber\\
&  \leq C\left\Vert \Lambda^{\frac{3}{2}}\theta(t)\right\Vert _{L^{2}%
}^{2-\frac{2\gamma}{9}}\left\Vert \theta_{0}\right\Vert _{L^{2}}^{\frac
{\gamma(9+2\gamma)}{9(3+\gamma)}}\left\Vert \Lambda^{\frac{3+\gamma}{2}}%
\theta(t)\right\Vert _{L^{2}}^{\frac{9+2\gamma}{3(3+\gamma)}}\nonumber\\
&  \leq C\left\Vert \Lambda^{\frac{3}{2}}\theta(t)\right\Vert _{L^{2}%
}^{4-\frac{4\gamma(6+\gamma)}{3(9+4\gamma)}}\left\Vert \theta_{0}\right\Vert
_{L^{2}}^{\frac{2\gamma(9+2\gamma)}{3(9+4\gamma)}}+\frac{1}{2}\left\Vert
\Lambda^{\frac{3+\gamma}{2}}\theta(t)\right\Vert _{L^{2}}^{2} \label{local6}%
\end{align}
which, in particular, gives
\begin{equation}
\left\Vert \theta(t)\right\Vert _{\dot{H}^{\frac{3}{2}}}\leq\frac{\left\Vert
\theta_{0}\right\Vert _{\dot{H}^{\frac{3}{2}}}}{\bigg[1-C\Big(2-\frac
{4\gamma(6+\gamma)}{3(9+4\gamma)}\Big)\left\Vert \theta_{0}\right\Vert
_{L^{2}}^{\frac{2\gamma(9+2\gamma)}{3(9+4\gamma)}}\left\Vert \theta
_{0}\right\Vert _{\dot{H}^{\frac{3}{2}}}^{2-\frac{4\gamma(6+\gamma
)}{3(9+4\gamma)}}t\bigg]^{\frac{1}{2-\frac{4\gamma(6+\gamma)}{3(9+4\gamma)}}}%
}, \label{local7}%
\end{equation}
for all $0\leq t\leq\bigg[2C\Big(2-\frac{4\gamma(6+\gamma)}{3(9+4\gamma
)}\Big)\left\Vert \theta_{0}\right\Vert _{L^{2}}^{\frac{2\gamma(9+2\gamma
)}{3(9+4\gamma)}}\left\Vert \theta_{0}\right\Vert _{\dot{H}^{\frac{3}{2}}%
}^{2-\frac{4\gamma(6+\gamma)}{3(9+4\gamma)}}\bigg]^{-1}$. \textit{A priori
}estimate (\ref{local7}) together with (\ref{max-L2}) yield
\[
\left\Vert \theta(t)\right\Vert _{H^{\frac{3}{2}}}\leq4\left\Vert \theta
_{0}\right\Vert _{H^{\frac{3}{2}}},\text{ for all }0\leq t\leq
\bigg[4C\left\Vert \theta_{0}\right\Vert _{L^{2}}^{\frac{2\gamma(9+2\gamma
)}{3(9+4\gamma)}}\left\Vert \theta_{0}\right\Vert _{\dot{H}^{\frac{3}{2}}%
}^{2-\frac{4\gamma(6+\gamma)}{3(9+4\gamma)}}\bigg]^{-1},
\]
and then $\left\Vert \theta(t)\right\Vert _{H^{\frac{3}{2}}}$ does not blow up
until
\begin{equation}
T_{1}=\frac{1}{C_{1}\left\Vert \theta_{0}\right\Vert _{L^{2}}^{\frac
{2\gamma(9+2\gamma)}{3(9+4\gamma)}}\left\Vert \theta_{0}\right\Vert _{\dot
{H}^{\frac{3}{2}}}^{2-\frac{4\gamma(6+\gamma)}{3(9+4\gamma)}}},
\label{Time-T1}%
\end{equation}
where $C_{1}$ is a constant that can be bounded from below and above for all
$\gamma\in\lbrack\gamma_{0},1]$. Thus, in (\ref{Dong-local-2}) we can consider
$T_{1}$ as in (\ref{Time-T1}).

On the other hand, using Gagliardo-Nirenberg inequality, we can estimate
$T^{\ast}$ in (\ref{timeeventual}) as
\[
T^{\ast}=C\alpha^{\frac{1}{1-\gamma}}\left\Vert \theta_{0}\right\Vert
_{L^{\infty}(\mathbb{T})}^{\frac{\gamma}{1-\gamma}}\leq t_{\gamma},
\]
where
\[
t_{\gamma}=Ck_{0}{}^{\frac{\gamma}{1-\gamma}}\alpha^{\frac{1}{1-\gamma}%
}\left\Vert \theta_{0}\right\Vert _{\dot{H}^{\frac{3}{2}}}^{\frac{\gamma
}{3(1-\gamma)}}\left\Vert \theta_{0}\right\Vert _{L^{2}}^{\frac{2\gamma
}{3(1-\gamma)}}%
\]
with $k_{0}$ independent of $\gamma,\alpha$ and $\theta_{0}$.

We claim that for $\gamma$ sufficiently close to $1$ we can choose $\alpha
\in(1-\gamma,1/2]$ such that $t_{\gamma}<T_{1}.$ Taking $C_{2}=C_{1}$ $k_{0}%
{}^{\frac{\gamma}{1-\gamma}}$, this is equivalent to
\begin{equation}
C_{0}C_{2}\alpha^{\frac{1}{1-\gamma}}\left\Vert \theta_{0}\right\Vert
_{\dot{H}^{\frac{3}{2}}}^{2-\frac{4\gamma(6+\gamma)}{3(9+4\gamma)}%
+\frac{\gamma}{3(1-\gamma)}}\left\Vert \theta_{0}\right\Vert _{L^{2}}%
^{\frac{2\gamma(9+2\gamma)}{3(9+4\gamma)}+\frac{2\gamma}{3(1-\gamma)}}\leq1.
\label{rcondition-0}%
\end{equation}
\bigskip Next, let us choose $R>0$ sufficiently large so that
\begin{equation}
\left\Vert \theta_{0}\right\Vert _{\dot{H}^{\frac{3}{2}}(\mathbb{T})}%
^{1-\frac{2\gamma}{3}}\left\Vert \theta_{0}\right\Vert _{L^{2}(\mathbb{T}%
)}^{\frac{2\gamma}{3}}\leq\left\Vert \theta_{0}\right\Vert _{H^{\frac{3}{2}%
}(\mathbb{T})}\leq R. \label{rcondition}%
\end{equation}
In view of (\ref{rcondition-0}) and (\ref{rcondition}), taking $C_{3}%
=C_{0}C_{2}$ it is sufficient to have that
\[
C_{3}\alpha^{\frac{1}{1-\gamma}}R^{\frac{18-3\gamma-2\gamma^{2}}%
{(1-\gamma)(9+4\gamma)}}\leq1
\]
or, equivalently,
\begin{equation}
\alpha\leq R^{-{\frac{18-3\gamma-2\gamma^{2}}{9+4\gamma}}}C_{3}^{-(1-\gamma)}.
\label{estimativefinal2}%
\end{equation}
Choosing $\alpha=\min\left\{  2(1-\gamma),\frac{1}{2}\right\}  $, it follows
from (\ref{estimativefinal2}) that there exists $\gamma_{1}:=\gamma_{1}%
(R)\in\lbrack\gamma_{0},1)$ such that $T^{\ast}\leq t_{\gamma}<T_{1}$ for all
$\gamma\in\lbrack\gamma_{1},1)$, which gives the claim.

Next, let $T_{\max}$ be the maximal existence time for the solution
(\ref{Dong-local-2}) of (\ref{1Dmodel}). Assume by contradiction that
$T_{\max}<\infty$. We have that $\theta\in C([0,T_{\max});H^{\frac{3}{2}%
}(\mathbb{T})))$ with $T^{\ast}<T_{1}\leq T_{\max}$. Then, by Theorem
\ref{eventual}, $\theta\in C^{\infty}(\mathbb{T\times}(T^{\ast},T_{\max}])$
and, in particular, $\theta(T_{\max})\in H^{\frac{3}{2}}(\mathbb{T}).$ So, by
using standard arguments and the local-existence of \cite{Dong}, we can extend
$\theta$ in the class (\ref{Dong-local-2}) to a time-interval $[0,T_{2})$ with
$T_{\max}<T_{2},$ which is a contradiction. It follows that $T_{\max}=\infty$
and $\theta$ is a global $H^{\frac{3}{2}}$-solution (which is classical) for
(\ref{1Dmodel}).

\fin


\begin{thebibliography}{99}                                                                                               %


\bibitem {Baker-Morlet-1}G. R. Baker, X. Li and A. C. Morlet, Analytic
structure of 1D-transport equations with nonlocal fluxes, Physica D 91 (1996), 349--375.

\bibitem {Bae}H. Bae, R. Granero-Belinch\'{o}n and O. Lazar, On the local and
global existence of solutions to 1D transport equations with nonlocal
velocity, arXiv:1806.01011, 2018.

\bibitem {Bae2}H. Bae, R. Granero-Belinch\'{o}n and Omar Lazar, Global
existence of weak solutions to dissipative transport equations with nonlocal
velocity, Nonlinearity 31 (2018), 1484--1515.

\bibitem {Calderon}A. Calder\'{o}n and A. Zygmund, Singular integrals and
periodic functions, Studia Math. 14 (1954), 249--271.

\bibitem {Constantin1}P. Constantin, A. Tarfulea and V. Vicol, Long time
dynamics of forced critical SQGV, Commun. Math. Phys. 335 (2015), 93--141.

\bibitem {Constantin2}P. Constantin and J. Wu, Regularity of H\"{o}lder
continuous solutions of the supercritical quasi-geostrophic equation, Annales
de l'Institut Henri Poincare Non Linear Analysis 25 (6) (2008), 1103--1110.

\bibitem {Cordoba}A. C\'{o}rdoba, D. C\'{o}rdoba and M. Fontelos, Formation of
Singularities for a Transport Equation with Nonlocal Velocity, Annals of
Mathematics 162 (3) (2005) 1377--1389.

\bibitem {Cordoba2}A. C\'{o}rdoba and D. C\'{o}rdoba, A maximum principle
applied to quasi-geostrophic equations, Commun. Math. Phys. (2004) (249) 511--528.

\bibitem {Vicol}M. Coti Zelati and V. Vicol, On the global regularity for the
supercritical SQG equation. Indiana Univ. Math. J. 65 (2) (2016), 535--552.

\bibitem {Do}T. Do, On a 1D transport equation with nonlocal velocity and
supercritical dissipation, Journal of Differential Equations, 256 (9) (2014), 3166--3178.

\bibitem {Dong}H. Dong, Well-posedness for a transport equation with nonlocal
velocity, Journal of Functional Analysis 255 (11) (2008), 3070--3097.

\bibitem {Ju}N. Ju, Dissipative 2D quasi-geostrophic equation: local
well-posedness, global regularity and similarity solutions, Indiana University
mathematics journal, (2007), 187--206.

\bibitem {Kiselev}A. Kiselev, Nonlocal maximum principles for active scalars,
Advances in Mathematics, 227 (5) 2011, 1806--1826.

\bibitem {Kiselev2}A. Kiselev, Regularity and blow up for active scalars.
Math. Model. Nat. Phenom. 5 (4) (2010), 225--255.

\bibitem {Li}D. Li and J. Rodrigo, Blow-up of solutions for a 1D transport
equation with nonlocal velocity and supercritical dissipation, Advances in
Mathematics, 217 (6) (2008), 2563--2568.

\bibitem {Silvestre-Vicol}L. Silvestre and V. Vicol, On a transport equation
with nonlocal drift, Trans. Amer. Math. Soc. 368 (2016), 6159--6188
\end{thebibliography}
\end{document}